
\def\LaTeX{%
  \let\Begin\begin
  \let\End\end
  \let\salta\relax
  \let\finqui\relax
  \let\futuro\relax}

\def\UK{\def\our{our}\let\sz s\def\analogue{analogue}}
\def\USA{\def\our{or}\let\sz z\def\analogue{analog}}



\LaTeX

\USA


\salta

\documentclass[twoside,a4paper,12pt]{article}
\setlength{\textheight}{24cm}
\setlength{\textwidth}{16cm}
\setlength{\oddsidemargin}{2mm}
\setlength{\evensidemargin}{2mm}
\setlength{\topmargin}{-15mm}
\parskip2mm


\usepackage[usenames,dvipsnames]{color}
\usepackage{amsmath}
\usepackage{amsthm}
\usepackage{amssymb}
\usepackage[mathcal]{euscript}
%
%
\usepackage{cite}
%
%
%


\definecolor{viola}{rgb}{0.3,0,0.7}
\definecolor{ciclamino}{rgb}{1,0,1}
\definecolor{rosso}{rgb}{0.8,0,0}

\def\gianni #1{{\color{blue}#1}}
\def\pier #1{{\color{rosso}#1}}

\def\bettibis #1{{\color{red}#1}}
\def\revis #1{{\color{red}#1}}

\def\gianni #1{#1}
\def\pier #1{#1}

\def\bettibis #1{#1}
\def\revis #1{#1}


\bibliographystyle{plain}


%

\finqui

\def\Beq{\Begin{equation}}
\def\Eeq{\End{equation}}
\def\Bsist{\Begin{eqnarray}}
\def\Esist{\End{eqnarray}}

\def\Bthm{\Begin{theorem}}
\def\Ethm{\End{theorem}}

\def\Brem{\Begin{remark}\rm}
\def\Erem{\End{remark}}

\def\Bdim{\Begin{proof}}
\def\Edim{\End{proof}}
\def\Bcenter{\Begin{center}}
\def\Ecenter{\End{center}}
\let\non\nonumber




\def\step #1 \par{\medskip\noindent{\bf #1.}\quad}


\def\Lip{Lip\-schitz}
\def\Holder{H\"older}

\def\aand{\quad\hbox{and}\quad}

\def\lhs{left-hand side}
\def\rhs{right-hand side}
\def\sfw{straightforward}


\def\organiz{organi\sz}


\def\multibold #1{\def\arg{#1}%
  \ifx\arg\pto \let\next\relax
  \else
  \def\next{\expandafter
    \def\csname #1#1#1\endcsname{{\bf #1}}%
    \multibold}%
  \fi \next}

\def\pto{.}

\def\multical #1{\def\arg{#1}%
  \ifx\arg\pto \let\next\relax
  \else
  \def\next{\expandafter
    \def\csname cal#1\endcsname{{\cal #1}}%
    \multical}%
  \fi \next}


\def\multimathop #1 {\def\arg{#1}%
  \ifx\arg\pto \let\next\relax
  \else
  \def\next{\expandafter
    \def\csname #1\endcsname{\mathop{\rm #1}\nolimits}%
    \multimathop}%
  \fi \next}

\multibold
qwertyuiopasdfghjklzxcvbnmQWERTYUIOPASDFGHJKLZXCVBNM.

\multical
QWERTYUIOPASDFGHJKLZXCVBNM.

\multimathop
dist div dom meas sign Sign supp .

\def\VV{{\mathbb{V}}}
\def\HH{{\mathbb{H}}}
\def\AA{{\mathbb{A}}}

\def\II{{\mathbb{I}}}
\def\FF{{\mathbb{F}}}


\def\accorpa #1#2{\eqref{#1}--\eqref{#2}}
\def\Accorpa #1#2 #3 {\gdef #1{\eqref{#2}--\eqref{#3}}%
  \wlog{}\wlog{\string #1 -> #2 - #3}\wlog{}}


\def\separa{\noalign{\allowbreak}}

\def\sign{\mathop{\rm sign}\nolimits}

\def\graffe #1{\mathopen\{#1\mathclose\}}

\def\<#1>{\mathopen\langle #1\mathclose\rangle}
\def\norma #1{\mathopen \| #1\mathclose \|}

\def\iot {\int_0^t}

\def\intQt{\int_{Q_t}}

\def\iO{\int_\Omega}

\def\dt{\partial_t}
\def\dn{\partial_n}

\def\cpto{\,\cdot\,}

\def\checkmmode #1{\relax\ifmmode\hbox{#1}\else{#1}\fi}
\def\aeO{\checkmmode{a.e.\ in~$\Omega$}}
\def\aeQ{\checkmmode{a.e.\ in~$Q$}}

\def\aet{\checkmmode{a.e.\ in~$(0,T)$}}

\def\aat{\checkmmode{for a.a.~$t\in(0,T)$}}


\def\erre{{\mathbb{R}}}




\def\genspazio #1#2#3#4#5{#1^{#2}(#5,#4;#3)}
\def\spazio #1#2#3{\genspazio {#1}{#2}{#3}T0}

\def\L {\spazio L}
\def\H {\spazio H}
\def\W {\spazio W}

\def\C #1#2{C^{#1}([0,T];#2)}


\def\Lx #1{L^{#1}(\Omega)}
\def\Hx #1{H^{#1}(\Omega)}

\def\LxG #1{L^{#1}(\Gamma)}
\def\HxG #1{H^{#1}(\Gamma)}

\def\Cx #1{C^{#1}(\overline\Omega)}

\def\LQ #1{L^{#1}(Q)}

\def\Luno{\Lx 1}
\def\Ldue{\Lx 2}
\def\Linfty{\Lx\infty}

\def\Huno{\Hx 1}
\def\Hdue{\Hx 2}
\def\Hunoz{{H^1_0(\Omega)}}
\def\Hmu{\Hx{-1}}

\def\Vz{V_0}
\def\Wz{W_0}


\def\LQ #1{L^{#1}(Q)}


\let\theta\vartheta
\let\eps\varepsilon
\let\phi\varphi

\let\TeXchi\chi                         
\newbox\chibox
\setbox0 \hbox{\mathsurround0pt $\TeXchi$}
\setbox\chibox \hbox{\raise\dp0 \box 0 }
\def\chi{\copy\chibox}



\def\phiz{\phi_0}

\def\phistar{\phi^*}

\def\rhostar{\rho^*}
\def\xistar{\xi^*}

\def\Omegaz{\Omega_0}

\def\normaH #1{\norma{#1}_H}

\let\hat\widehat
\def\Beta{\hat{\vphantom t\smash\beta\mskip2mu}\mskip-1mu}
\def\Pi{\hat\pi}
\def\betaz{\beta^\circ}
\def\betaeps{\beta_\eps}
\def\signeps{\sign_\eps}

\def\Betaeps{\Beta_\eps}
\def\phieps{\phi_\eps}
\def\mueps{\mu_\eps}
\def\sigmaeps{\sigma_\eps}
\def\xieps{\xi_\eps}
\def\zetaeps{\zeta_\eps}

\def\Ieps{I_\eps}
\def\hphi{\hat\phi}
\def\hsigma{\hat\sigma}
\def\veps{v_\eps}

\def\Tstar{T^*}


\def\g #1{\gamma_{#1}}
\def\muG{\mu_\Gamma}
\def\muH{\mu_{\calH}}
\def\sigmas{\sigma_{\!s}}
\def\sigmaz{\sigma_0}
\def\sigmastar{\sigma^*}
\def\soluz{(\phi,\mu,\sigma,\xi,\zeta)}
\def\soluzeps{(\phieps,\mueps,\sigmaeps)}

\def\Csh{C_{sh}}
\def\Csys{C_{sys}}
\def\hatC{\hat C}

\def\Mz{M_0}
\def\Mpi{M_\pi^*}

\def\Arho{A(\rho)}

\def\Lpi{L_\pi}

\def\vz{v_0}

\Begin{document}


\title{Sliding mode control for \revis{a} phase field system\\
related to tumor growth
}

\author{}
\date{}
\maketitle
\Bcenter
\vskip-2cm
{\large\sc Pierluigi Colli$^{(1)}$}\\
{\normalsize e-mail: {\tt pierluigi.colli@unipv.it}}\\[.25cm]
{\large\sc Gianni Gilardi$^{(1)}$}\\
{\normalsize e-mail: {\tt gianni.gilardi@unipv.it}}\\[.25cm]
{\large\sc Gabriela Marinoschi$^{(2)}$}\\
{\normalsize e-mail: {\tt gabriela.marinoschi@acad.ro}}\\[.25cm]
{\large\sc Elisabetta Rocca$^{(1)}$}\\
{\normalsize e-mail: {\tt elisabetta.rocca@unipv.it}}\\[.45cm]
$^{(1)}$
{\small Dipartimento di Matematica ``F. Casorati'', Universit\`a di Pavia}\\
{\small and \ IMATI-C.N.R., Pavia}\\
{\small via Ferrata 1, 27100 Pavia, Italy}\\[.2cm]
$^{(2)}$
{\small ``Gheorghe Mihoc-Caius Iacob'' Institute of Mathematical Statistics\\
and \ Applied Mathematics of the Romanian Academy}\\
{\small Calea 13 Septembrie 13, 050711 Bucharest, Romania}\\[.2cm]
\Ecenter

\Begin{abstract}
In the present contribution we study the sliding mode control (SMC) problem for a diffuse interface tumor growth model coupling a viscous Cahn--Hilliard type equation for the phase variable with a reaction-diffusion equation for the nutrient. 
First, we prove the well-posedness and some regularity results for the state system 
modified by the state-feedback control law. 
Then, we show that the chosen SMC law forces the system to reach within finite time the sliding manifold
(that we chose in order that the tumor phase remains constant in time). 
The feedback control law is added in the Cahn--Hilliard type equation and leads  the phase  
onto a prescribed target~$\phi^*$ in finite time. 
\vskip3mm

\noindent {\bf Key words:} sliding mode control, 
Cahn--Hilliard system, reaction-diffusion equation,  tumor growth, nonlinear boundary value problem, state-feedback control law.
\vskip3mm
\noindent {\bf AMS (MOS) Subject Classification:} 34H15, 35K25, 35K61, 93B52, 92C50, 97M60.
\End{abstract}

\salta

\pagestyle{myheadings}
\newcommand\testopari{\sc Colli \ --- \ Gilardi \ --- \ Marinoschi \ --- \ Rocca}
\newcommand\testodispari{\sc Sliding modes for a tumor growth model}
\markboth{\testodispari}{\testopari}

\finqui


\section{Introduction}
\label{Intro}
\setcounter{equation}{0}

Sliding mode control (SMC) - which is today considered a classic instrument for regulation
of continuous or discrete systems in finite-dimensional settings
(see e.g.~the monographs \cite{BFPU08, EFF06, ES99, FMI11, I76, Utkin92, UGS09, YO99}) - has been acknowledged  as one of the
basic approaches to the design of robust controllers for nonlinear complex
dynamics that work under uncertainty. One of the main example of  complex systems studied nowadays
both in biomedical and mathematical literatures is related to the tumor growth dynamics. For the case of an incipient tumor, i.e., before the development of quiescent cells, the studied diffuse interface type phase field models often consist of a Cahn--Hilliard equation coupled with a reaction-diffusion equation for the nutrient (cf., e.g., \cite{CLLW09, GLSS16, HZO12, HKNZ15}). In this work, in particular,  we consider the problem of sliding mode control for a tumor growth model recently introduced in \cite{GLSS16}. In comparison with \cite{GLSS16}, we have neglected here the effects of chemotaxis and active transport, but the new feature of \eqref{Isecondabis} is the inclusion of the SMC law $\rho \sign(\phi-\phistar)$, where $\rho$ is a positive parameter that will be chosen large enough. This term  forces the system trajectories onto
the sliding surface  $\phi=\phistar$ in finite time. All in all, we consider here  the following viscous Cahn--Hilliard/Reaction-Diffusion model for tumor growth
\Bsist
  & \dt\phi - \Delta\mu = (\g1 \sigma - \g2) p(\phi)
  & \quad \hbox{in $Q:=\Omega\times(0,T)$}
  \label{Iprima}
  \\
  &  \mu = \tau \dt\phi - \Delta\phi + F'(\phi) + \rho \sign(\phi-\phistar)
& \quad \hbox{in $Q$}
  \label{Isecondabis}  
  \\
  & \dt\sigma - \Delta\sigma = - \g3 \sigma p(\phi) + \g4 (\sigmas - \sigma) + g 
  & \quad \hbox{in $Q$}
  \label{Iterza}
\Esist
where $\Omega$ is the domain in which the evolution takes place,
$T$~is some final time,
$\phi$ denotes the difference in volume fraction, where $\varphi = 1$ represents the tumor phase and $\varphi = -1$ represents the healthy tissue phase, $\mu$~is the chemical potential
and $\sigma$ is the concentration of a nutrient for the tumor cells (e.g., oxigen or glucose).
Moreover, $\tau$~is a positive viscosity coefficient,
$\g i$~for $i=1,\dots,4$ denotes the positive constant proliferation rate, apoptosis rate, nutrient consumption rate, and nutrient supply rate, respectively. The term $ \g 1 p(\phi)\sigma$ models the proliferation of tumor cells which is proportional to the concentration of the nutrient, the term $ \g 2p(\phi)$ models the apoptosis of tumor cells, and $\g 3 p(\phi) \sigma$ models the consumption of the nutrient only by the tumor cells. The constant $\sigmas$ denotes the nutrient concentration in a pre-existing vasculature, and $\g 4(\sigmas - \sigma)$ models the supply of nutrient from the blood vessels if $\sigmas > \sigma$ and the transport of nutrient away from the domain $\Omega$ if $\sigmas< \sigma$. The function $g$ is a source term which may represent the supply of a nutrient (see \cite{BC}), or even a drug in a  chemotherapy.
Moreover, 
$F'$~stands for the derivative of a double-well potential~$F$
and $p$ is a {smooth} nonnegative proliferation function on the domain of~$F$.
Typical examples of potentials, meaningful in view of applications, are
\begin{align}
  & F_{reg}(r) = \frac 14 \, (r^2-1)^2\,,
  \quad r \in \erre
  \label{regpot}
  \\
  & F_{log}(r) = (1+r)\ln (1+r)+(1-r)\ln (1-r) - c_0 \, r^2 \,,
  \quad r \in (-1,1)
  \label{logpot}
  \\
  & F_{obs}(r) = I(r) - c_0 \, r^2 \,,
  \quad r \in \erre
  \label{obspot}
\end{align}
where $c_0>1$ in \eqref{logpot} in order to produce a double-well,
while $c_0$ is an arbitrary positive number in~\eqref{obspot},
and the function $I$ in~\eqref{obspot} is the indicator function of~$[-1,1]$, i.e.,
it takes the values $0$ or $+\infty$ according to whether or not $r$ belongs to~$[-1,1]$.
The {potentials} \eqref{regpot} and \eqref{logpot} are 
the classical regular potential and the so-called logarithmic potential, respectively.
More generally, the potential $F$ could be just the sum $F=\Beta+\Pi$,
where $\Beta$ is a convex function that is allowed to take the value~$+\infty$,
and $\Pi$ is a smooth perturbation (not necessarily concave).
In such a case, $\Beta$~is supposed to be proper and lower semicontinuous
so that its subdifferential is well defined while  the derivative
might not exist.
This happens, for example,  in the case~\eqref{obspot}
and then equation \eqref{Isecondabis} becomes a differential inclusion.
Finally, the operator $\sign:\erre\to 2^{\erre}$ is 
defined~by
\Beq
  \sign r := \frac r {|r|}
  \quad \hbox{if $r\not=0$}
  \aand
  \sign 0 := [-1,1] .
  \label{defsignintro}
\Eeq
The aim  of introducing such a feedback law in \eqref{Isecondabis}  is to force the order parameter to reach a prescribed distribution $\phistar$ in a finite time.
However, the resulting problem of forcing the solution of the modified system
to reach the manifold $\phi=\phistar$ in a finite time looks difficult.
Indeed, we can ensure the existence of the desired sliding mode only under 
a suitable compatibility condition between the measure of the set~$\Omega$ and the viscosity 
parameter $\tau$ (cf.~the following \eqref{Ismallness}).

The above system is complemented by initial conditions like $\phi(0)=\phiz$ and $\sigma(0)=\sigmaz$
and suitable boundary conditions.
Concerning the latter, 
we take the usual homogeneous Neumann conditions for~$\phi$ and~$\sigma$, that~is,
\Beq
  \dn\phi = 0
  \aand
  \dn\sigma = 0
  \quad \hbox{on $\Sigma := \Gamma\times (0,T) $}
  \non
\Eeq
where $\Gamma$ is the boundary of~$\Omega$
and $\dn$ is the (say, outward) normal derivative.
Instead, we consider a Dirichlet boundary condition for the chemical potential,~i.e.,
\Beq
  \mu = \muG
  \quad \hbox{on $\Sigma$}
  \label{Imubou}
\Eeq
where $\muG$ is a given smooth function.
This choice is twofold: from one side it looks reasonable from the modelling point of view , in case $\muG=0$,  the condition \eqref{Imubou} allows for the free flow of
cells across the outer boundary (cf.~\cite{WLFC08} and \cite{BCT} where similar conditions are imposed on a chemical potential in a different framework). On the other hand,  we also need \eqref{Imubou} for the analysis. Indeed, a major technical issue here, in case we choose a usual Neumann boundary condition for $\mu$, would be to estimate its mean value. This would be doable if $\rho=0$ and the potential $F$ is assumed to have a controlled growth (cf., e.g., \cite{GLR17}), but it is not the case when the feedback law is added in \eqref{Isecondabis} and under our fully general assumptions on $F$. 

At this point, without the aim of completeness, let us describe the recent literature both on tumor growth modelling and on SMC related to our problem. 

Modelling tumor growth dynamics has recently become a major issue in applied mathematics (see, e.g.,  \cite{CL10,WLFC08}). Numerical simulations of diffuse interface models for tumor growth have been carried out in several papers (see, e.g., \cite[Ch. 8]{CL10}); nonetheless, a rigorous mathematical analysis of the resulting systems of PDEs is still in its infancy (cf., e.g., \cite{CGH, CGRS1, CGRS2, DFRSS17, FGR, FLR, GLDarcy, GLDirichlet, GLNeumann}). 
Recently,  in \cite{WLFC08} the authors introduced a continuum diffuse interface model of multispecies tumor growth and tumor-induced angiogenesis in two and three dimensions for investigating their morphological evolution. They make use of the Cahn--Hilliard framework which originated from the theory of phase transitions, and which is used extensively in materials science and multiphase fluid flow.
Other diffuse interface models including chemotaxis and transport effects  have been subsequently introduced (cf.~\cite{GLNS, GLSS16}) and also the formal sharp interface limits have been investigated. Rigorous sharp interface limits have been obtained in some particular cases in the two recent works \cite{MR, RS}. 

Regarding the SMC  literature, the SMC scheme is well known for its robustness against variations of dynamics, disturbances,
time-delays and nonlinearities. The design procedure
of a SMC system is a two-stage process. The first phase is to choose a set of sliding manifolds
such that the original system restricted to the intersection of \gianni{them} has a desired
behavior. In this paper,  we choose to force the tumor phase parameter to stay constant in time within finite time with the obvious application in mind that the phase $\phi$ should become as close as possible to the constant value $\phi=-1$ corresponding to the case when no tumorous phase is present anymore or to a configuration $\phistar$ which is suitable for surgery.  The second step is to design a SMC law that forces the system trajectories to stay onto
the sliding surface. To this end, we have added the term $ \rho \sign(\phi-\phistar)$ in the Cahn--Hilliard evolution for $\phi$ (cf.~\eqref{Isecondabis}) in order to force $\phi$ to stay equal to a given desired value $\phistar$ in a finite time. 

Sliding mode controls are pretty attractive in many applications.
As a result, in recent years there has been a growing interest in extending well-developed methods
for finite-dimensional systems described by ODEs
(see, for example, \cite{LO02, O83, O00, OU83})
or also to control infinite-dimensional dynamical systems (cf.~\cite{OU83, OU87, OU98}). 
Moreover, the theoretical development in Hilbert spaces or for PDE systems has only taken the attention in the last ten years.
In this regard,  we can cite the papers \cite{CRS11, Levaggi13, XLGK13}
concerning the control of semilinear PDE systems.

Finally,  we can quote the recent contribution \cite{BCGMR}, where a sliding mode approach is applied for the first time to phase field systems of Caginalp type coupling the evolution of a phase paratemer to the one of the relative temperature,  and the chosen SMC laws force the system to reach within finite time the sliding manifold. In that case it was possible to have different choices for the manifold:  in particular, either one of the physical variables or a combination of them could remain constant in time. With reference to the results of \cite{BCGMR}, we aim to mention the analyses developed in \cite{Colt1, Colt2}: in particular, the second contribution is devoted to a conserved phase field system with a SMC feedback law for the internal energy in the temperature equation. 

In the present contribution, instead, we are forced,  mainly from the fact that we have a fourth order equation for $\phi$, to include a sliding mode control of the type $\rho \sign(\phi-\phistar)$ in the chemical potential $\mu$ (cf.~\eqref{Isecondabis}) and we cannot handle neither the presence of the $\sigma$-dependence in it or a non-local in space control law (as we did in some cases in \cite{BCGMR}). We need here, for technical reasons, to include a local in space control, i.e., such that
its value at any point and any time just depends on the value of the state. In this way, however,  as already mentioned above, we need to enforce a compatibility condition bewteen the size of $|\Omega|$ and the viscosity coefficient $\tau$ of the type 
\Beq
2 \Csh \,  |\Omega|^{2/3}
  < \tau
  \label{Ismallness}
\Eeq
where $\Csh$ is the constant related to some embedding inequalities (cf.~\eqref{smallness} in the next Section~\ref{STATEMENT}). 
Such a condition means that either $|\Omega|$ has to be sufficiently small
once the shape of $\Omega$ is fixed in the sense of the following Remark~\ref{Cshape} or $\tau$ must be sufficiently large compared to the size of $|\Omega|$. Regarding the fact that we do not treat here a feedback law depending also on $\sigma$, this turns out to be quite reasonable in view of applications since we mainly aim to optimize the tumor cell distributions and not the nutrient concentration, in general. 

Other approaches to the problem of control for tumor growth models are possible, even if a few mathematical results are presently available on this subject in the literature. In the two recent papers \cite{GLR17} and \cite{CGRS17} the authors face the problem of finding  first order necessary optimality conditions for the minimization of a cost functional forcing the phase to approach the desired target $\phistar$ in the best possible way by means of a control  variable representing the concentration of cytotoxic drugs in \cite{GLR17} and the supply of a nutrient or a drug in a  chemotherapy, in \cite{CGRS17} (it could be the function $g$ in \eqref{Iterza} in the present contribution). 

The main advantage of controlling the sliding mode is that it strengthens
the trajectories of the system to reach the sliding surface and keep it on \gianni{it in} a pointwise way, while, in general, within the classical optimal control theory (cf., e.g., \cite{GLR17, CGRS17}), one can get just necessary optimality conditions and the control is nonlocal in space and/or in time. 

The paper is \organiz ed as follows.
In the next section, we list our assumptions, state the problem in a precise form
and present our results.
The last two sections are devoted to the corresponding proofs.
Section~\ref{WELLPOSEDNESS} deals with well-posedness,
while the existence of the sliding {mode} is proved in Section~\ref{SLIDINGMODES}.


\section{Statement of the problem and results}
\label{STATEMENT}
\setcounter{equation}{0}

In this section, we describe the problem under study
and present our results.
As in the Introduction,
$\Omega$~is the body where the evolution takes place.
We assume $\Omega\subset\erre^3$
to~be open, bounded, connected, and smooth,
and we write $|\Omega|$ for its Lebesgue measure.
Moreover, $\Gamma$ and $\dn$ still stand for
the boundary of~$\Omega$ and the outward normal derivative, respectively.
Given a finite final time~$T>0$,
we set for convenience
$Q:=\Omega\times(0,T)$.
Furthermore, if $X$ is a Banach \bettibis{space}, the symbol $\norma\cpto_X$ denotes its norm,
with the exception of the norms in the $L^\infty$ spaces on $\Omega$, $Q$ and~$\Sigma$,
for which we use the same symbol $\norma\cpto_\infty$ since no confusion can arise.
Finally, the dual space of $X$ is denoted by~$X^*$
and we write $\<\cpto,\cpto>_{\pier{{X^*, X}}}$ for the duality paring between $X^*$ and~$X$.

Now, we specify the assumptions on the structure of our system.
We assume that
\Bsist
  && \g i \in [0,+\infty) \quad \hbox{for $i=1,2,3$}, \quad
  \g4 ,\, \tau \in (0,+\infty)
  \aand
  \sigmas \in \erre
  \label{hpconst}
  \\
  && \Beta : \erre \to [0,+\infty]
  \quad \hbox{is convex, proper and l.s.c.}
  \quad \hbox{with} \quad
  \Beta(0) = 0
  \label{hpBeta}
  \\
  && \Pi: \erre \to \erre
  \quad \hbox{is a $C^1$ function}
  \aand
  {\Pi\,}'
  \quad \hbox{is \Lip\ continuous},
  \label{hpPi}
\\
&& \pier{ p: \erre \to [0,+\infty)
  \quad \hbox{is a bounded and Lipschitz continuous  function}.}
  \label{hp-p}
\Esist
\Accorpa\HPstruttura hpconst hp-p
We set for brevity
\Beq
  \beta := \partial\Beta , \quad
 \pier{ \pi := {\Pi\,}',
  \quad
  \Lpi = \hbox{the \Lip\ constant of $\pi$}}
  \label{defbetapi}
\Eeq
and denote by $D(\beta)$ and $D(\Beta)$
the effective domains of $\beta$ and~$\Beta$, respectively.
Next, in order to simplify notations, we~set
\Bsist
  && H := \Ldue, \quad
  V := \Huno, \quad
  \Vz := \Hunoz
  \label{defVH}
  \\
  && W := \graffe{v\in\Hx2: \dn v=0}
  \aand
  \Wz := \Hx2 \cap \Hunoz
  \label{defW}
\Esist
\Accorpa\Defspazi defVH defW
and endow these spaces with their standard norms.
Moreover, we denote by $\Csh$ a constant realizing the inequalities
\Bsist
  && \norma v_\infty
  \leq \Csh \, |\Omega|^{1/6} \normaH{\Delta v}
  \quad \hbox{for every $v\in\Wz$}
  \label{embDir}
  \\
  && \norma v_\infty
  \leq \Csh \bigl( |\Omega|^{-1/2} \normaH v + |\Omega|^{1/6} \normaH{\Delta v} \bigr)
  \quad \hbox{for every $v\in W$} \,.
  \label{embNeu}
\Esist

\Brem
\label{Cshape}
We show that such a constant actually exists 
and depends on $\Omega$ just through its shape. 
\revis{%
Hence, we consider a class of open sets having the same shape.
To do this, we fix an open set $\Omegaz\subset\erre^3$ with Lebesgue measure~$1$.
Then, each open set of the same class is related to $\Omegaz$ by the formula
$\Omega=\{x_0\}+\lambda R\,\Omegaz$, where $x_0$ is a point in $\erre^3$, 
the real number $\lambda$ is positive and $R$ belongs to the rotation group $SO(3)$. 
Let now $\Csh$ be a constant satisfying \accorpa{embDir}{embNeu}
with $\Omega$ and $|\Omega|$ replaced by $\Omegaz$ and~$1$, respectively.
Such a constant exists since the three-dimensional open set $\Omegaz$ 
is supposed to be bounded and smooth, as usual.
Now, we take any $v\in W$ and check~\eqref{embNeu}.
We define $\vz$ belonging to the analogue of $W$ constructed on~$\Omegaz$, i.e.,
$\vz\in H^2(\Omegaz)$ with $\dn\vz=0$, by the formula
$\vz(y):=v(x_0+\lambda Ry)$ for $y\in\Omegaz$.
Then, it is \sfw\ to show that
\Bsist
  & \norma\vz_{L^\infty(\Omegaz)}
  = \norma v_{\Linfty} 
  \aand
  \norma\vz_{L^2(\Omegaz)}
  = \lambda^{-3/2} \norma v_{\Ldue}
  \label{normevz}
  \\
  & \norma{\Delta\vz}_{L^2(\Omegaz)}
  = \lambda^{1/2} \norma{\Delta v}_{\Ldue} .
  \label{normaDvz}
\Esist
On the other hand, we have that $\lambda=|\Omega|^{1/3}$
(choose $v\equiv1$, whence $\vz\equiv1$, in the second~\eqref{normevz}).
Therefore, we deduce that
\Bsist
  && \norma v_{\Linfty}
  = \norma\vz_{L^\infty(\Omegaz)}
  \leq \Csh \bigl( \norma\vz_{L^2(\Omegaz)} + \norma{\Delta\vz}_{L^2(\Omegaz)} \bigr)
  \non
  \\
  && = \Csh \bigl( |\Omega|^{-1/2} \norma v_{\Ldue} + |\Omega|^{1/6} \norma{\Delta v}_{\Ldue} \bigr)
  \non
\Esist
i.e., \eqref{embNeu}.
The derivation of \eqref{embDir} is similar and even simpler.}
\Erem

At this point, we describe the state system 
modified by the state-feedback control law.
We introduce the operator $\sign:\erre\to 2^{\erre}$
defined~by
\Beq
  \sign r := \frac r {|r|}
  \quad \hbox{if $r\not=0$}
  \aand
  \sign 0 := [-1,1] .
  \label{defsign}
\Eeq
Notice that $\sign$ is the subdifferential of the real function
$r\mapsto|r|$ and thus is maximal monotone.
Next, we reduce the Dirichlet boundary condition $\mu=\muG$ to the \bettibis{homogeneous} one.
By assuming $\muG\in\L2{\HxG{1/2}}$ just to start with,
we introduce the harmonic extension $\muH$ of $\muG$ defined \pier{in this way}:
\aat, $\muH(t)$ is the unique solution to the problem
\Beq
  \muH(t) \in \Huno , \quad
  - \Delta\muH(t) = 0 \quad \hbox{in \pier{$\calD' (Q)$}}
  \aand
  \muH(t)|_\Gamma = \muG(t) \,.
  \label{defmuH}
\Eeq
Then, we take $\mu-\muH$ as new \bettibis{unknown}.
However, in order to avoid a new notation, we still term $\mu$ the above difference.
Thus, the problem to be solved is the following: \pier{we are}
given the functions $g$, $\muG$, $\phistar$ and the initial data $\phiz$ and $\sigmaz$ such that
\Bsist
  && g \in \LQ\infty , \quad
  \muG \in \H1{\LxG2} \cap \L\infty{\HxG{3/2}} 
  \label{hpdati}
  \\
  && \phistar \in W
  \aand
  \inf D(\beta) < \inf\phistar \leq \sup\phistar < \sup D(\beta)
  \label{hpstar}
  \\
  && \phiz \in W , \quad
  \betaz(\phiz) \in H
  \aand
  \sigmaz \in V \cap \Linfty
  \label{hpz}
\Esist
\Accorpa\Hpdati hpdati hpz
\pier{where $\betaz$ denotes the minimal section of $\beta.$} \bettibis{Notice that the assumptions on $\muG$ in \eqref{hpdati} are additional and do not follow from \revis{what is previously stated} on $\muH$}. Then,
we look for a quintuplet $\soluz$ satisfying the regularity requirements
\Bsist
  && \phi \in \W{1,\infty}H \cap \H1V \cap \L\infty W
  \label{regphi}
  \\
  && \mu \in \L\infty\Wz
  \label{regmu}
  \\
  && \sigma \in \H1H \cap \L\infty V \cap \L2W 
  \label{regsigma}
  \\
  && \xi \in \L\infty H
  \aand
  \zeta \in \L\infty H
  \label{regxizeta}
\Esist
\Accorpa\Regsoluz regphi regxizeta
and solving
\Bsist
  && \dt\phi - \Delta\mu = (\g1 \sigma - \g2) \, p(\phi)
  \quad \aeQ
  \label{prima}
  \\
  && \mu = \tau \dt\phi - \Delta\phi + \xi + \pi(\phi) + \rho \, \zeta - \muH
  \quad \aeQ
  \label{seconda}
  \\
  && \dt\sigma - \Delta\sigma = - \g3 \sigma \, p(\phi) + \g4 (\sigmas - \sigma) + g 
  \quad \aeQ
  \label{terza}
  \\
  && \xi \in \beta(\phi) 
  \aand
 \pier{\zeta} \in \sign(\phi - \phistar)
  \quad \aeQ
  \label{quarta}
  \\
  && \phi(0) = \phiz 
  \aand
  \sigma(0) = \sigmaz
  \label{cauchy}
\Esist
\Accorpa\Pbl prima cauchy
where $\rho$ is a positive parameter.
We notice that the boundary conditions $\dn\phi=0$, $\mu=0$ and $\dn\sigma=0$
are contained in \eqref{regphi}, \eqref{regmu} and \eqref{regsigma}, respectively,
due to the definitions \eqref{defW} of $W$ and~$\Wz$.
We also remark that
\Beq
  \muG \in L^\infty(\Sigma)
  \label{muGbdd}
\Eeq
as a \pier{consequence} of \eqref{hpdati}.
Indeed, $\HxG{3/2}\subset\LxG\infty$ since $\Sigma$ is a two-dimensional smooth surface.
Here is our well-posedness result.

\Bthm
\label{Wellposedness}
Assume \HPstruttura\ and \Hpdati.
Then, for every $\rho>0$, there exists at least one quintuplet $\soluz$
fulfilling \Regsoluz, solving problem \Pbl\
and satisfying the estimates
\Bsist
  && \norma\mu_\infty
  \leq \Csh \, \frac {2 |\Omega|^{2/3}} \tau \, \rho
  + \hatC 
  \label{stimamu}
  \\
  && |\sigma| \leq \sigmastar := \max\{\norma{\sigmas+\g4^{-1}g}_\infty,\norma\sigmaz_\infty\}
  \quad \aeQ
  \label{maxprinc}
\Esist
where $\Csh$ is the same as in \accorpa{embDir}{embNeu}
and the constant $\hatC$ depends only on $\Omega$, $T$ 
and the quantities involved in assumptions \HPstruttura\ and \Hpdati.
In particular, $\hatC$~does not depend on~$\rho$.
Moreover, the components $\phi$ and $\sigma$ of the solution are uniquely determined.
\Ethm

The above result is quite general.
In particular, all the potentials \accorpa{regpot}{obspot} 
are certainly allowed.

As far as the problem of a sliding mode is concerned,
we prove a result that only involves the component $\phi$ of the solution,
which is uniquely determined.
However, we can ensure the existence of a sliding mode at least for $\rho$ large enough
only under a restriction.
Namely, we need the following condition
\Beq
  \Csys := 
  \Csh \, \frac {2 |\Omega|^{2/3}} \tau
  < 1 
  \label{smallness}
\Eeq
where $\Csh$ is the constant that appears in \accorpa{embDir}{embNeu}. 
Such a condition means that $|\Omega|$ has to be sufficiently small
once the shape of $\Omega$ is fixed in the sense of Remark~\ref{Cshape}.

\Bthm
\label{Sliding}
In addition to \HPstruttura\ and \Hpdati,
assume that
\Beq
  \Delta\phistar \in \Linfty \pier{.}
  \label{hpsliding}
\Eeq
Moreover, assume~\eqref{smallness}.
Then, there exists $\rhostar>0$, 
depending only on $\Omega$, $T$, the structure and the data of the problem,
such that, for every $\rho>\rhostar$, the following is true:
if $\soluz$ is a solution to problem \Pbl\ whose component $\mu$ satisfies~\eqref{stimamu},
there exists a time $\Tstar\in[0,T)$ such~that
\Beq
  \phi(t) = \phistar
  \quad \aeO
  \quad \hbox{for every $t\in[\Tstar,T]$}.
  \label{sliding}
\Eeq
In particular, there exists a solution for which \eqref{sliding} holds true.
\Ethm

\Brem
\label{Calcolo}
In the proof we give in Section~\ref{SLIDINGMODES},
we show that possible values of $\rhostar$ and $\Tstar$ that fit the above statement~are
\Bsist
  && \rhostar =
  \frac 1 {1-\Csys} \Bigl(
    \hatC 
    + M
    + \Mpi
    + \frac \tau T \, \Mz  
  \Bigr)
  \aand
  \Tstar = \frac \tau {\rho-\Arho} \, \Mz
  \non
\Esist
where $M$, $\Mz$, $\Mpi$ and $\Arho$ are given by
\Bsist
  && M :=
  \norma\muG_\infty
  + \norma{\Delta\phistar}_\infty
  + \norma\xistar_\infty \,, \quad
  \Mz := \norma{\phiz-\phistar}_\infty 
  \non
  \\
  && \Mpi := \sup\{|\pi(\phistar(x)+r)| :\ x\in\Omega,\ |r|\leq\Mz\}
  \non
  \\
  && A(\rho) := 
  \Csys \, \rho
  + \hatC
  + M 
  + \Mpi
  \quad \hbox{for $\rho>0$}
  \non
\Esist
\pier{where $\xistar := \betaz(\phistar)$.} In these formulas, $\hatC$ is the same as in \eqref{stimamu}.
We will see that the above definitions ensure that $\Arho<\rho$ for $\rho>\rhostar$
and that $\Tstar\in[0,T)$.
\Erem

The rest of the section is devoted to make some notations precise
and to introduce some tools we use in the remainder of the paper.
In performing our a priori estimates, 
we often account for the \Holder\ inequality and the Young inequality
\Beq
  ab \leq \delta a^2 + \frac 1 {4\delta} \, b^2
  \quad \hbox{for every $a,b\geq0$ and $\delta>0$}.
  \label{young}
\Eeq
Moreover, we repeatedly use the notation
\Beq
  Q_t := \Omega \times (0,t)
  \aand
  \Sigma_t := \Gamma \times (0,t)
  \quad \hbox{for $t\in(0,T)$}.
  \label{defQ}
\Eeq
For simplicity, we usually omit $dx$, $ds$, etc.\ in integrals.
More precisely, we explicitly write, e.g., $ds$ only if the variable $s$ actually appears 
in the function under the integral sign.
We also take advantage of the Dirichlet problem solver operator $\calD:\Hmu\to\Hunoz$
defined as follows:
if $f\in\Hmu$, then $\calD f$ is the unique solution $u$ to the Dirichlet problem
\Beq
  u \in \Hunoz
  \aand
  - \Delta u = f \,.
  \label{defD}
\Eeq
As $\Omega$ is bounded and smooth, we have that
$\calD f\in\Wz$ as well~as
\Beq
  \norma{\calD f}_{\Hdue} \leq C \normaH f
\pier{, \quad \hbox{in particular,}\quad }
  \normaH{\calD f} \leq C \normaH f
  \quad \hbox{for every $f\in H$}
  \label{regD}
\Eeq
where $C$ depends only on~$\Omega$.
Furthermore, we define the equivalent norm in~$\Hmu$ by the formula
\Beq
  \norma f_*^2 := \iO |\nabla\calD f|^2
  \quad \hbox{for $f\in\Hmu$}.
  \label{defnormastar}
\Eeq
Notice that
\Beq
  \< f , \calD f >_{\pier{\Vz^*, \Vz}} = \norma f_*^2
 \aand
 \pier{ \< f(t) , \calD (\dt f(t)) >_{\pier{\Vz^*, \Vz}}
  = \frac 12 \, \frac d{dt} \, \norma{f(t)}_*^2, \quad t\in (0,T), }
  \label{propnormastar}
\Eeq
for every $f\in\Hx{-1}$ and $f\in\H1{\Hx{-1}}$, respectively.
For the same reason as above, the assumptions on $\muG$ in \eqref{hpdati} imply
a proper regularity for $\muH$ and the corresponding estimates.
Namely, \pier{owing to the maximum principle and to the elliptic regularity results stated, e.g., in \cite[Thm.~3.2, p.~1.79]{BG87}}, we have~that
\Bsist
  && \norma\muH_\infty
  \leq \norma\muG_\infty
  \label{maxmuH}
  \\
  && \norma\muH_{\L\infty\Hdue} 
  \leq C \norma\muG_{\L\infty{\HxG{3/2}}} 
  \label{stimamuH}
  \\
  && \norma{\dt\muH}_{\L2H}
  \leq C \norma{\dt\muG}_{\L2{\LxG2}} 
  \label{stimadtmuH}
\Esist
with a constant $C$ that depends only on~$\Omega$. 

Finally, while a particular care is taken in computing some constants,
we follow a general rule to denote less important ones,
in order to avoid boring calculations.
The small-case symbol $c$ stands for different constants which depend only
on~$\Omega$, on the final time~$T$, the shape of the nonlinearities
and on the constants and the norms of
the functions involved in the assumptions of our statements, but~$\rho$.
The dependence on $\rho$ will be always written explicitly, indeed.
Constants depending on further parameters are characterized by a corresponding subscript.
So, e.g., we write $c_\eps$ for constants depending on $\eps$ as well.
Hence, the meaning of $c$ and $c_\eps$ might change from line to line 
and even in the same chain of equalities or inequalities.
On the contrary, we mark precise constants which we can refer~to
by using different symbols, e.g., capital letters,
mainly with indices.


\section{Well-posedness}
\label{WELLPOSEDNESS}
\setcounter{equation}{0}

This section is devoted to the proof of Theorem~\ref{Wellposedness}.
We first prove the partial uniqueness given in the statement.
Then, we introduce a proper regularization and construct a solution
satisfying both \eqref{stimamu} and the maximum principle~\eqref{maxprinc}.
It is convenient to observe once and for all that every solution satisfies
\Beq
  \calD(\dt\phi)
  + \tau \dt\phi
  - \Delta\phi
  + \xi
  + \pi(\phi)
  + \rho \, \zeta
  - \muH
  = \calD \bigl( (\g1 \sigma - \g2) \, p(\phi) \bigr)
  \quad \aeQ \,.
  \label{nomu}
\Eeq
Indeed, it sufficies to apply $\calD$ to both sides of \eqref{prima}
(by~recalling that $\mu$ is $\Vz$-valued)
and replace $\mu$ by means of~\eqref{seconda}.

\step
Partial uniqueness

We pick two solutions $(\phi_i,\mu_i,\sigma_i,\xi_i,\zeta_i)$, $i=1,2$,
\pier{corresponding to the same data}
and prove that $(\phi_1,\sigma_1)=(\phi_2,\sigma_2)$.
We write \eqref{nomu} and \eqref{terza} for both solutions and take the difference.
By setting for brevity $\phi:=\phi_1-\phi_2$ 
and defining the analogous differences $\mu$, $\sigma$, $\xi$ and~$\zeta$,
we have that \aeQ
\Bsist
  && \calD(\dt\phi)
  + \tau \dt\phi
  - \Delta\phi
  + \xi
  + \rho \, \zeta
  \non
  \\
  && = \pi(\phi_2) - \pi(\phi_1)
  + \calD \bigl( \g1 \sigma \, p(\phi_1) \bigr)
  + \calD \bigl( (\g1 \sigma_2 - \g2) (p(\phi_1) - p(\phi_2)) \bigl)
  \non
  \\
  && \dt\sigma - \Delta\sigma
  + \g4 \sigma
  = - \g3 \sigma \, p(\phi_1)
  - \g3 \sigma_2 \bigl( p(\phi_1) - p(\phi_2) \bigr).
  \non
\Esist
Now, we multiply such equations by $\phi$ and~$\sigma$, respectively,
sum up and integrate over~$Q_t$, where $t\in(0,T)$ is arbitrary.
\pier{Thanks} to the boundary and initial conditions, 
the definition \eqref{defnormastar} of~$\norma\cpto_*$,
the second inequality in~\eqref{regD},
the \Lip\ continuity of $\pi$ and~$p$,
the boundedness of~$p$,
and the inequality \eqref{maxprinc} satisfied by $\sigma_1$ and~$\sigma_2$,
we obtain
\Bsist
  && \frac 12 \, \norma{\phi(t)}_*^2
  + \frac \tau 2 \, \normaH{\phi(t)}^2
  + \intQt |\nabla\phi|^2
  + \intQt \xi \phi
  + \rho \intQt \zeta \phi
  \non
  \\
  && \quad {}
  + \frac 12 \, \normaH{\sigma(t)}^2 
  + \intQt |\nabla\sigma|^2
  + \g4 \intQt |\sigma|^2
  \non
  \\
  && \leq c \intQt \bigl( |\phi|^2 + |\phi| \, |\sigma| + |\sigma|^2 \bigr)
  \leq c \intQt \bigl( |\phi|^2 + |\sigma|^2 \bigr)
  \non
\Esist
where the values of $c$ might depend on the solutions we are considering as well.
All the terms on the \lhs\ are nonnegative
(those involving $\xi$ and $\zeta$ by monotonicity).
Therefore, by applying the Gronwall lemma,
we conclude that $\phi=\sigma=0$, i.e., that $(\phi_1,\sigma_1)=(\phi_2,\sigma_2)$.

\medskip

At this point, we start proving the existence of a solution 
satisfying the properties of the statements.
To this end, we introduce the approximating problem
obtained by replacing the graphs $\sign$ and $\beta$
by the \Lip\ continuous functions $\signeps$ and~$\betaeps$, 
their Yosida regularizations at the level~$\eps\in(0,1)$
(\pier{concerning general properties of maximal monotone and 
subdifferential operators along with their Yosida regularizations, 
the reader can see, e.g., \cite{Barbu,Brezis}}).
We also introduce their primitives $|\cdot|_\eps$ and~$\Betaeps$ vanishing at the origin.
Such functions satisfy the properties listed below
(the last one also being a consequence of the assumptions~\eqref{hpBeta})\pier{:}
\Bsist
  && \signeps r = \frac r {\max\{\eps,|r|\}}
  \quad \hbox{for every $r\in\erre$}
  \label{formulasigneps}
  \\
  && 0 \leq |r|_\eps := \int_0^r \signeps s \, ds \leq |r|
  \quad \hbox{for every $r\in\erre$}
  \label{defmoduloeps}
  \\
  && |\betaeps(r)| \leq |\betaz(r)|
  \quad \hbox{for every $r\in D(\beta)$}, \quad \hbox{where}
  \non
  \\
  && \quad \betaz(r)\ \hbox{is the element of $\beta(r)$ having minimum modulus}
  \label{disugbetaeps}
  \\
  && 0 \leq \Betaeps(r) := \int_0^r \betaeps(s) \, ds \leq \Beta(r)
  \quad \hbox{for every $\, r\in D(\Beta)$}.
  \label{propBetaeps}
\Esist
Moreover, we choose the following regularization $\Ieps:\erre\to\erre$ of the identity:
\Beq
  \Ieps(r) := \pier{\max\{-1/\eps, \min \{r, 1/\eps\}\}}
  \quad \hbox{for $r\in\erre$}.
  \label{defIeps}
\Eeq
Thus, the approximating problem consists in finding a triplet $\soluzeps$ 
satisfying \Regsoluz\ and solving
\Bsist
  && \dt\phieps - \Delta\mueps = (\g1 \Ieps(\sigmaeps) - \g2) \, p(\phieps)
  \quad \aeQ
  \label{primaeps}
  \\
  && \mueps = \tau \dt\phieps - \Delta\phieps + \xieps + \pi(\phieps) + \rho \, \zetaeps - \muH
  \quad \aeQ
  \label{secondaeps}
  \\
  && \dt\sigmaeps - \Delta\sigmaeps = - \g3 \Ieps(\sigmaeps) \, p(\phieps) + \g4 (\sigmas - \sigmaeps) + g 
  \quad \aeQ
  \label{terzaeps}
  \\
  && \quad \hbox{where} \quad
  \xieps := \betaeps(\pier{\phieps}) 
  \aand
  \pier{\zetaeps} := \signeps(\phieps - \phistar)
  \label{quartaeps}
  \\
  && \phieps(0) = \phiz 
  \aand
  \sigmaeps(0) = \sigmaz \pier{.}
  \label{cauchyeps}
\Esist
\Accorpa\Pbleps primaeps cauchyeps
\pier{Also in this case, the boundary conditions are contained in the regularity requirements.}

\Bthm
\label{Wellposednesseps}
Under the assumptions of Theorem~\ref{Wellposedness},
the problem \Pbleps\ has a unique solution $\soluzeps$ 
satisfying the regularity \bettibis{requirements} \Regsoluz.
Moreover, the component $\sigmaeps$ fulfils the inequality
\Beq
  |\sigmaeps| \leq \sigmastar
  \quad \aeQ
  \label{maxprinceps}
\Eeq
with the same $\sigmastar$ as in~\eqref{maxprinc}.
\Ethm

\Bdim
We show the well-posedness of the approximating problem
in a less regular framework.
Namely, we look for solutions satisfying
\Bsist
  && \phieps \in \H1H \cap \L\infty V \cap \L2W
  \label{regphieps}
  \\
  && \mueps \in \L2\Wz
  \label{regmueps}
  \\
  && \sigmaeps \in \H1H \cap \L\infty V \cap \L2W 
  \label{regsigmaeps}
\Esist
instead of \Regsoluz.
In the sequel of the paper, it will be clear that the solution we find here
also satisfies the regularity we have required, 
even though we will proceed formally.
We present the approximating problem in the equivalent form obtained by replacing \eqref{primaeps}
by the analogue of~\eqref{nomu},~i.e.,
\Beq
  \calD(\dt\phieps)
  + \tau \dt\phieps
  - \Delta\phieps
  + \xieps
  + \pi(\phieps)
  + \rho \, \zetaeps
  - \muH
  = \calD \bigl( (\g1 \Ieps(\sigmaeps) - \g2) \, p(\phieps) \bigr)
  \label{nomueps}
\Eeq
where $\xieps$ and $\zetaeps$ are still given by~\eqref{quartaeps}.
We want to show that we can solve \eqref{nomueps}, \eqref{terzaeps} and \eqref{cauchyeps} for $(\phieps,\sigmaeps)$
with the proper regularity that is needed.
Then, we will use \eqref{secondaeps} as a definition of~$\mueps$.
In order to solve the sub-problem,
we present it as a Cauchy problem for a nonlinear abstract equation of the type
\Beq
  \frac d {dt} \, (\phieps,\sigmaeps)
  + \AA (\phieps,\sigmaeps)
  + \FF (\phieps,\sigmaeps)
  = G
  \label{astratta}
\Eeq
in the framework of the Hilbert triplet
$(\VV,\HH,\VV^*)$ where
\Beq
  \VV := V \times V
  \aand
  \HH := H \times H
  \label{hilbert}
\Eeq
with a non-standard embedding $\HH\subset\VV^*$.
To justify what we are going to write,
in particular the forthcoming definitions \eqref{defprod} and~\accorpa{defFuno}{defFdue},
we first notice that $\calD$ is a symmetric linear continuous operator from $H$ into itself
which satisfies the first \pier{property in}~\eqref{propnormastar}.
This implies that the integral $\iO(\calD v)v$ is nonnegative for every $v\in H$,
so that the operator $\calD+\tau I$, where $I$ is the identity map of~$H$, 
is an isomorphism from $H$ onto itself and $(\calD+\tau I)^{-1}$ 
is a well-defined linear continuous operator from $H$ into itself.
Once this is established,
we observe that the system given by \eqref{nomueps} and \accorpa{terzaeps}{quartaeps}
is equivalent to the following variational equation
\Bsist
  && \iO \bigl( \calD(\dt\phieps) + \tau\dt\phieps \bigr) v
  + \iO \dt\sigmaeps \, z
  + \iO \nabla\phieps \cdot \nabla v
  + \iO \nabla\sigmaeps \cdot \nabla z
  \non
  \\
  && \quad {}
  + \iO (\calD+\tau I) (\calD+\tau I)^{-1} \Bigl(
    (\betaeps+\pi+\rho\signeps)(\phieps)
    - \calD \bigl( (\g1 \Ieps(\sigmaeps) - \g2) \, p(\phieps) \bigr)
  \Bigr) v
  \non
  \\
  && \quad {}
  + \iO \bigl(
    \g3 \Ieps(\sigmaeps) \, p(\phieps)
    - \g4 (\sigmas - \sigmaeps)
  \bigr) z
  \non
  \\
  && = \iO \bigl( (\calD+\tau I) (\calD+\tau I)^{-1} \muH \bigr) \, v
  + \iO g z
  \non
\Esist
holding \aet, for every $(v,z)\in V\times V$.
This equation can be written as
\Bsist
  && \bigl( \dt(\phieps,\sigmaeps) , (v,z) \bigr)_{\HH} 
  + \< \AA (\phieps,\sigmaeps) , (v,z) >_{\pier{\VV^*, \VV}}
  + \bigl( \FF(\phieps,\sigmaeps),(v,z) \bigr)_{\HH}
  \non
  \\
  && = \bigl( ((\calD+\tau I)^{-1} \muH, g) , (v,z) \bigr)_{\HH}
  \non
\Esist
that is, in the form \eqref{astratta}
with an obvious meaning of~$G$,
provided that $(\cpto,\cpto)_{\HH}$ and $\AA:\VV\to\VV^*$ are given by
\Bsist
  && \bigl( (\hphi,\hsigma) , (v,z) \bigr)_{\HH}
  := \iO \bigl( \calD\hphi + \tau\hphi \bigr) v
  + \iO \hsigma z
  \quad \hbox{for $(\hphi,\hsigma),(v,z)\in\HH$}
  \label{defprod}
  \\
  && \< \AA (\hphi,\hsigma) , (v,z) >_{\pier{\VV^*, \VV}}
  := \iO \nabla\hphi \cdot \nabla v
  + \iO \nabla\hsigma \cdot \nabla z
  \quad \hbox{for $(\hphi,\hsigma),(v,z)\in\VV$}
  \qquad
  \label{defA}
\Esist
and the function $\FF:\HH\to\HH$ is defined by the following rule:
for $(\hphi,\hsigma)\in\HH$, the value
$\FF(\hphi,\hsigma)$ is the pair $\bigl(\FF_1(\hphi,\hsigma),\FF_2(\hphi,\hsigma)\bigr)$
given~by
\Bsist
  && \FF_1(\hphi,\hsigma)
  := (\calD+\tau I)^{-1} \bigl(
      (\betaeps+\pi+\rho\signeps)(\hphi)
      - \calD \bigl( (\g1 \Ieps(\hsigma) - \g2) \, p(\hphi) \bigr)
    \bigl)
  \qquad
  \label{defFuno}
  \\
  && \FF_2(\hphi,\hsigma)
  := \g3 \Ieps(\hsigma) \, p(\hphi)
    - \g4 (\sigmas - \hsigma).
  \label{defFdue}
\Esist
By owing to our prelimirary observations,
we see that, from one side, formula \eqref{defprod} actually yields an inner product in $\HH$ 
that is equivalent to the standard one;
on the other side, by also noting that the real functions
\Beq
  r \mapsto (\betaeps+\pi+\rho\signeps)(r),
  \quad r \in \erre ,
  \aand
  (r,s) \mapsto \Ieps(s) \, p(r),
  \quad (r,s) \in \erre^2
  \non
\Eeq
are \Lip\ continuous \pier{(cf. \eqref{hp-p})},
we conclude that $\FF:\HH\to\HH$ is \Lip\ continuous.
\revis{Note that, in particular, the boundedness of $p$ is used in order to check the \Lip\ continuity of the considered function of $(r,s)$.
In addition, we} observe at once that the meaning of the embedding $\II:\HH\to\VV^*$ is the following:
for arbitrary $(\hphi,\hsigma)\in\HH$ and $(v,z)\in\VV$, we have that
\Beq
  \< (\hphi,\hsigma) , (v,z) >_{\pier{\VV^*, \VV}}
  = \bigl( (\hphi,\hsigma) , (v,z) \bigr)_{\HH}
  = \iO \bigl( \calD\hphi + \tau\hphi \bigr) v
  + \iO \hsigma z .
  \non
\Eeq
Now, we show that $\AA:\VV\to\VV^*$ is maximal monotone.
Indeed, it is monotone since
\Beq
  \< \AA (v,z) , (v,z) >_{\pier{\VV^*, \VV}}
  = \iO (|\nabla v|^2 + |\nabla z|^2)
  \geq 0
  \quad \hbox{for every $(v,z)\in\VV$}.
  \non
\Eeq
On the other hand, $\AA+\II:\VV\to\VV^*$
is surjective since it holds for every $(v,z)\in\VV$
\Bsist
  && \< (\AA+\II) (v,z) , (v,z) >_{\pier{\VV^*, \VV}}
  = \< \AA(v,z) , (v,z) >_{\pier{\VV^*, \VV}}
  + \bigl( (v,z) , (v,z) \bigr)_{\HH}
  \non
  \\
  && = \iO (|\nabla v|^2 + |\nabla z|^2)
  + \iO \bigl( \calD v + \tau v \bigr) v
  + \iO |z|^2
  \geq \min\{1,\tau\} \norma{(v,z)}_{\VV}^2 \,.
  \non
\Esist
We deduce that also the realization of $\AA$ in~$\HH$
(or~restriction of $\AA$ to~$\HH$, which we still term~$\AA$) defined by setting
$D(\AA):=W\times W$ is maximal monotone.
The monotonicity property obviously follows.
To see the maximal monotonicity, we show that, for any $G=(G_1,G_2)\in\HH$,
we can find $(\hphi,\hsigma)\in W\times W$
satisfying $(\AA+\II)(\hphi,\hsigma)=G$.
Indeed, the equation to be solved has a solution $(\hphi,\hsigma)\in\VV$.
On the other hand, such a solution satisfies
\Beq
  \< \AA(\hphi,\hsigma) , (v,z) >_{\pier{\VV^*, \VV}}
  + \bigl( (\hphi,\hsigma) , (v,z) \bigr)_{\HH}
  = \bigl( (G_1,G_2) , (v,z) \bigr)_{\HH}
  \quad \hbox{for every $(v,z)\in\VV$}
  \non
\Eeq
that is the couple of variational Neumann problems
\Bsist
  && \iO \nabla\hphi \cdot \nabla v
  + \tau \iO \hphi \, v
  = \iO (\calD G_1 + \tau \revis{G_1} - \calD \hphi) v
  \quad \hbox{for every $v\in V$}
  \non
  \\
  && \iO \nabla\hsigma \cdot \nabla z
  + \iO \hsigma z
  = \iO G_2 z
  \quad \hbox{for every $z\in V$}
  \non
\Esist
so that both $\hphi$ and $\hsigma$ belong to $W$ by the elliptic regularity theory.
Once such a maximal monotonicity is established,
we recall that the Cauchy problem to be solved is the equation \eqref{astratta}
complemented by the initial condition~\eqref{cauchyeps}.
Thus, we can apply the general theory and ensure that such a problem 
has a unique solution satisfying~\eqref{regphieps} and \eqref{regsigmaeps}.
At this point, \eqref{secondaeps}~yields $\mueps$ with a level of regularity
that is lower than required.
However, from \eqref{nomueps}~and \pier{\eqref{secondaeps}}
we deduce that $\mueps$ solves \eqref{primaeps} as well,
so that also the regularity requirement \eqref{regmueps} holds for~$\mueps$.
This concludes the proof of well-posedness.

Finally, we prove~\eqref{maxprinceps}.
In order to show the upper inequality $\sigma\leq\sigmastar$,
we rewrite \eqref{terzaeps} in the form
\Bsist
  && \dt\sigmaeps
  - \Delta\sigmaeps
  + \g3 \, p(\phieps) \bigl( \Ieps(\sigmaeps) - \Ieps(\sigmastar) \bigr)
  + \g4 (\sigmaeps - \bettibis{\sigmastar}) 
  \non
  \\
  && = - \g3 \, p(\phieps) \, \Ieps(\sigmastar)
  + \bigr(\g4 \sigmas + g - \g4 \sigmastar \bigr)
  \non
\Esist
and observe that the \rhs\ is nonpositive
due to the definition of $\sigmastar$ contained in \eqref{maxprinc}
and the positivity assumptions on the constants and on~$p$.
Thus, multiplying by the positive part $\veps:=(\sigmaeps-\sigmastar)^+$,
integrating over~$Q_t$, accounting for the boundary and initial conditions,
and observing that $\veps(0)=0$,
we obtain for every $t\in[0,T]$
\Beq
  \frac 12 \iO |\veps(t)|^2
  + \intQt |\nabla\veps|^2
  + \g3 \intQt p(\phieps) \bigl( \Ieps(\sigmaeps) - \Ieps(\sigmastar) \bigr) (\sigmaeps-\sigmastar)^+
  + \g4 \intQt |\veps|^2
  \leq 0 \,.
  \non
\Eeq
By also recalling the definition \eqref{defIeps} of~$\Ieps$,
we see that every term on the \lhs\ is nonnegative.
Hence, we deduce that $\veps=0$, i.e., $\bettibis{\sigmaeps}\leq\sigmastar$.
In order to derive the inequality $\bettibis{\sigmaeps}\geq-\sigmastar$,
one writes \eqref{terza}~as
\Bsist
  && \dt\sigmaeps
  - \Delta\sigmaeps
  + \g3 \, p(\phieps) \bigl( \Ieps(\sigmaeps) + \Ieps(\sigmastar) \bigr)
  + \g4 (\sigmaeps + \bettibis{\sigmastar}) 
  \non
  \\
  && = \g3 \, p(\phieps) \, \Ieps(\sigmastar)
  + \bigr(\g4 \sigmas + g + \g4 \sigmastar \bigr)
  \non
\Esist
and multiplies by $-(\sigmaeps+\sigmastar)^-$.
\Edim

Our project is now the following:
we show that the approximating solution satisfies some bounds
(in particular, we formally derive the further regularity \Regsoluz\ for $\soluzeps$);
then, we take the limit as $\eps\searrow0$ by accounting for compactness and monotoniciy arguments.
As far as the notation for the \pier{constants} is concerned,
we follow the general rule explained at the end of the previous section.
In particular, the (possibly different) values denoted by $c$
are independent of $\rho$ and~$\eps$.
The same holds for the precise constants $C_1$, $C_2$, etc., 
that we introduce in the sequel and mark with capital letters for possible future references.
\pier{Moreover, due to \eqref{maxprinceps} and the definition \eqref{defIeps}, 
from now on in \eqref{primaeps}, \eqref{terzaeps},
\eqref{nomueps} we can replace $\Ieps(\sigmaeps)$ by $ \sigmaeps$,
 assuming that \gianni{$1/\eps \geq \sigmastar$}, that is, $ \eps \leq 1/\sigmastar$.} 
\bettibis{Finally, let us note that from now on we will underline the dependence in the estimates from $\tau$ and $|\Omega|$ only 
when it is necessary in order to conclude the proofs.} 

\step
First a priori estimate

\pier{In view of}~\eqref{maxprinceps}, 
we notice that the \rhs\ of \eqref{terzaeps} is uniformly bounded in $\LQ2$
by a constant depending only on our structural assumptions and on the norms
$\norma g_\infty$ and $\norma\sigmaz_\infty$.
By the parabolic regularity theory, we infer that
\Beq
  \norma\sigmaeps_{\H1H\cap\L\infty V\cap\L2W} \leq C_1 \,.
  \label{primastima}
\Eeq

\step
Second a priori estimate

We \pier{add $\phieps$ to both sides of \eqref{nomueps}, then we test by $\dt\phieps$ and integrate over~$(0,t)$}
by accounting for the boundary and initial conditions,
and rearrange.
By recalling the first property \pier{of the norm~$\norma\cpto_*$ in \eqref{propnormastar}
and the} formulas \eqref{propBetaeps} and~\eqref{formulasigneps},
we obtain
\Bsist
  && \iot \norma{\dt\phieps(s)}_*^2
  + \tau \intQt |\dt\phieps|^2 
  + \frac 12 \pier{\norma{\phieps(t)}_V^2}
  \non
  \\
  && \quad {}
  + \iO \Betaeps(\phieps(t))
  + \rho \iO |\pier{\phieps(t) - \phistar}|_\eps
  \non
  \\
  \separa
  && = \frac 12 \pier{\norma{\phiz}_V^2}
  + \iO \Betaeps(\phiz)
  + \rho \iO |\pier{\phiz - \phistar}|_\eps
  \non
  \\
  && \quad {}
  + \intQt \calD \bigl( (\g1 \sigmaeps - \g2) \, p(\phieps) \bigr) \, \dt\phieps
  + \intQt \muH \dt\phieps
  \pier{{}+ \intQt (\phieps - \pi(\phieps))} \dt\phieps \,.
  \non
\Esist
Every term on the \lhs\ is nonnegative.
\pier{Using} the Young inequality, the structural \pier{assumptions~\Hpdati}, \eqref{stimamuH}, and the estimates \eqref{regD} and~\eqref{maxprinceps},
we see that the \rhs\ can be bounded from above~by
\Bsist
  && \frac 12 \, \pier{\norma{\phiz}_V^2}
  + \iO \Beta(\phiz)
  + \rho \iO  |\pier{\phiz - \phistar}|
  + \frac \tau 2 \intQt |\pier{\dt\phieps}|^2 
  \non
  \\
  && \quad {}
  + \pier{c{}} \, \norma{\calD \bigl( (\g1 \sigmaeps - \g2) \, p(\phieps) \bigr) + \muH}_{\LQ2}^2
  + c \intQt (1 + |\phieps|^2)
  \non
  \\
  && \leq \frac \tau 2 \intQt |\pier{\dt\phieps}|^2 
  + c + c \intQt |\phieps|^2 . 
  \non
\Esist
\pier{Hence, by applying} the Gronwall lemma, we conclude that
\Beq
  \norma\phieps_{\H1H\cap\L\infty V}
  \leq C_2 (\rho^{1/2} + 1) .
  \label{secondastima}
\Eeq
We also get \pier{the same bound for $\Betaeps(\phieps)$ in $\L\infty\Luno$.}

\step
Third a priori estimate

We (formally) differentiate \eqref{nomueps} with respect to time
(by accounting for \eqref{quartaeps}).
We obtain
\Bsist
  && \calD(\dt^2\phieps)
  + \tau \dt^2\phieps
  - \Delta\dt\phieps
  \non
  \\
  && \quad {}
  + \betaeps'(\phieps) \dt\phieps
  + \pi'(\phieps) \dt\phieps
  + \rho \, \signeps'(\phieps-\phistar) \dt(\phieps-\phistar)
  - \dt\muH
  \non
  \\
  && = \calD \bigl( \pier{\g1 \dt\sigmaeps } \, p(\phieps) + (\g1 \sigmaeps - \g2) p'(\phieps) \dt\phieps \bigr).
  \non
\Esist
Now, we \pier{test} this equation by $\dt\phieps$ and \pier{integrate over~$(0,t)$}.
By recalling the second \pier{identity in} \eqref{propnormastar} and rearranging, we have that
\Bsist
  && \frac 12 \, \norma{\dt\phieps(t)}_*^2
  + \frac\tau 2 \iO |\dt\phieps(t)|^2
  + \intQt |\nabla\dt\phieps|^2
  \non
  \\
  && \quad {}
  + \intQt \betaeps'(\phieps) |\dt\phieps|^2
  + \rho \intQt \signeps'(\phieps-\phistar) |\dt(\phieps-\phistar)|^2
  \non
  \\
  && = \frac 12 \, \norma{\dt\phieps(0)}_*^2
  + \frac\tau 2 \iO |\dt\phieps(0)|^2
  \non
  \\
  && \quad {}
  + \intQt \bigl\{
    \calD \bigl( \pier{\g1 \dt\sigmaeps  \, p(\phieps)}+ (\g1 \sigmaeps - \g2) p'(\phieps) \dt\phieps \bigr)
  \bigr\} \dt\phieps
  \non
  \\
  && \quad {}
  - \intQt \pi'(\phieps) |\dt\phieps|^2
  + \intQt \dt\muH \, \dt\phieps \,.
  \non
\Esist
\pier{All the terms} on the \lhs\ are nonnegative.
Moreover, the second inequality in~\eqref{regD},
the estimates \eqref{primastima} and \eqref{secondastima},
and the assumption on $\muG$ in \eqref{hpdati} combined with \eqref{stimamuH}
allow to find a bound for the volume integrals on the \rhs.
Therefore, we deduce that
\Bsist
  && \frac\tau 2 \iO |\dt\phieps(t)|^2
  + \intQt |\nabla\dt\phieps|^2
  \non
  \\
  && \leq \frac 12 \, \norma{\dt\phieps(0)}_*^2
  + \frac\tau 2 \iO |\dt\phieps(0)|^2
  + c (\rho + 1)
  \label{perterzastima}
\Esist
and it remains to estimate the norms of $\dt\phieps(0)$.
To this end, we write \eqref{nomueps} at the time $t=0$,~i.e.,
\Bsist
  && \calD(\dt\phieps(0))
  + \tau \dt\phieps(0)
  \non
  \\
  && = \calD \bigl( (\g1 \sigmaz - \g2) \, p(\phiz) \bigr)
  + \Delta\phiz
  - \betaeps(\phiz)
  - \rho \, \signeps(\phiz)
  - \pi(\phiz)
  + \muH(0)
  \non
\Esist
and \pier{test} by $\dt\phieps(0)$.
\pier{By recalling} the first \pier{identity in} \eqref{propnormastar},
\eqref{regD}, \eqref{disugbetaeps} and our assumptions on the data,
we easily infer (with the help of the Young inequality) that
\Beq
  \norma{\dt\phieps(0)}_*^2
  + \tau \iO |\dt\phieps(0)|^2
  \leq 
  \frac\tau 2 \iO |\dt\phieps(0)|^2
  + \frac {|\Omega|} \tau \, \rho^2 + c \,.
  \non
\Eeq
By combining this with \eqref{perterzastima}, we deduce that
\Beq
  \frac\tau 2 \iO |\dt\phieps(t)|^2
  + \intQt |\nabla\dt\phieps|^2
  \leq \frac {|\Omega|} \tau \, \rho^2
  + c(\rho+1)
  \non
\Eeq
\aat\ and conclude that
\Beq
  \norma{\dt\phieps}_{\L\infty H}
  \leq \frac {(2|\Omega|)^{1/2}} \tau \, \rho
  + C_3 \bigl( \rho^{1/2} + 1 \bigr)
  \aand
  \norma{\dt\phieps}_{\L2V}
  \leq C_3' (\rho + 1).
  \label{terzastima}
\Eeq

\step
Fourth a priori estimate

We write \eqref{primaeps} in the form
\Beq
  - \Delta\mueps = (\g1 \sigmaeps - \g2) \, p(\phieps) - \dt\phieps 
  \non
\Eeq
and \gianni{owe to} \eqref{terzastima}, \eqref{maxprinceps} and the Young inequality 
\gianni{(using $2^{1/2} <2$)}
to deduce that
\Bsist
  && \norma{\Delta\mueps}_{\L\infty H} 
  \leq \norma{\dt\phieps}_{\L\infty H}
  + c (\norma\sigmaeps_\infty + 1)
  \non
  \\
  && \leq \frac {(2|\Omega|)^{1/2}} \tau \, \rho
  + c \bigl( \rho^{1/2} + 1 \bigr)
  \leq \frac {2 |\Omega|^{1/2}} \tau \, \rho 
  + c \,.
  \non
\Esist
In particular, we have that \pier{(cf.~\eqref{defD} and~\eqref{regD})}
\Beq
  \norma\mueps_{\L\infty\Hdue} \leq c(\rho+1).
  \label{stimamuepsHdue}
\Eeq
\pier{Moreover,} we \pier{state}  a more precise $L^\infty$~estimate
\pier{on account of}~\eqref{embDir}.
Namely, we have that
\Beq
  \norma\mueps_\infty
  \leq \Csh \, \frac {2 |\Omega|^{2/3}} \tau \, \rho
  + C_4 \,.
  \label{quartastima}
\Eeq

\step
Fifth a priori estimate

We write \eqref{secondaeps} in the form
\Beq
  - \Delta\phieps
  + \xieps
  + \rho \, \zetaeps 
  = \mueps
  - \tau \dt\phieps
  - \pi(\phieps)
  + \muH  
  \non
\Eeq
and multiply by~$-\Delta\phieps$.
In dealing with the third term on the \lhs\
we split $-\Delta\phieps$ as $-\Delta(\phieps-\phistar)-\Delta\phistar$.
By integrating over~$\Omega$ and rearranging, we \pier{infer that}
\Bsist
  && \iO |{-}\Delta\phieps|^2
  + \iO \betaeps'(\phieps) |\nabla\phieps|^2
  + \pier{\rho} \iO \signeps'(\phieps-\phistar) |\nabla(\phieps-\phistar)|^2
  \non
  \\
  && = \iO \bigl(
    \mueps
    - \tau \dt\phieps
    - \pi(\phieps)
    + \muH 
  \bigr)(-\Delta\phieps) 
  + \rho \iO \zetaeps \, \Delta\phistar 
  \non
\Esist
\pier{ \aet.} We recall that $\betaeps'$ and $\signeps'$ are nonnegative 
and that $|\zetaeps|\leq1$ \aeQ.
Thus, using the Young and Schwarz inequalities, we deduce that
\Beq
  \frac 12 \iO |{-}\Delta\phieps|^2
  \leq \frac 12 \, \normaH{\mueps-\tau\dt\phieps-\pi(\phieps)+\muH}^2
  + \rho \, |\Omega|^{1/2} \normaH{\Delta\phistar}
  \non
\Eeq
whence
\Bsist
  \normaH{-\Delta\phieps \pier{(t)} }
  &\leq& \normaH{(\mueps-\tau\dt\phieps-\pi(\phieps)+\muH  ) \pier{(t)}}
  + c \, \rho^{1/2}\non
\\
  &\leq& \normaH{\mueps  \pier{(t)}}
  + \tau \normaH{\dt\phieps \pier{(t)}}
  + c \,\pier{ ( \rho^{1/2} + 1) \quad \hbox{for a.e. } t\in (0,T)}
  \non
\Esist
thanks to \eqref{secondastima} and the \Lip\ continuity of~$\pi$.
\pier{Owing} to \eqref{quartastima} for~$\mueps$ and \eqref{terzastima} for~$\dt\phieps$,
we have that
\Beq
  \norma{-\Delta\phieps}_{\L\infty H}
  \leq c \, \rho
  + c (\rho^{1/2}+1)
  \leq c (\rho+1).
  \non
\Eeq
\pier{Hence, due to \eqref{secondastima} and to} the regularity theory for elliptic equations, we deduce that
\Beq
  \norma\phieps_{\L\infty W}
  \leq C_5 (\rho+1).
  \label{quintastima}
\Eeq
\bettibis{Finally, by comparison we infer that
\Beq
\norma\xieps_{\L\infty H}\leq C_6(\rho+1).
 \label{sestastima}
\Eeq
}

\step
Conclusion

At this point, we \pier{collect} all the previous estimates 
and \pier{recall} that $|\zetaeps|\leq1$ \aeQ.
Then, using standard compactness results \pier{and possibly taking a subsequence},
we obtain
\Bsist
  & \phieps \to \phi
  & \quad \hbox{weakly star in $\W{1,\infty}H\cap\H1V\cap\L\infty W$}
  \qquad
  \label{convphi}
  \\
  & \mueps \to \mu
  & \quad \hbox{weakly star in $\L\infty\Wz$}
  \label{convmu}
  \\
  & \sigmaeps \to \sigma
  & \quad \hbox{weakly star in $\H1H\cap\L\infty V\cap\L2W$}
  \label{convsigma}
  \\
  & \xieps \to \xi
  & \quad \hbox{weakly star in $\L\infty H$}
  \label{convxi}
  \\
  & \zetaeps \to \zeta
  & \quad \hbox{weakly star in $\LQ\infty$}
  \label{convzeta}
\Esist
for some quintuplet $\soluz$ satisfying \Regsoluz.
Moreover, the estimates we have obtained are conserved for the limiting functions.
In particular, those given by 
\eqref{quartastima} and \eqref{maxprinceps}
are conserved, that is, the triplet $(\phi,\mu,\sigma)$ satisfies
\accorpa{stimamu}{maxprinc} with $\hatC=C_4$.
Furthermore, \pier{in the light of} the compact embeddings $\Hdue\subset V \pier{{}\subset H}$ and $\Hdue\subset\Cx0$
and taking advantage, e.g., of \cite[Sect.~8, Cor.~4]{Simon},
we derive strong convergence for $\phieps$ and~$\sigmaeps$.
Namely, we have that
\Bsist
  & \phieps \to \phi 
  & \quad \hbox{strongly in $\C0V\cap C^0(\overline Q)$}
  \label{strongphi}
  \\
  & \sigmaeps \to \sigma 
  & \quad \hbox{strongly in \pier{$\C0H \cap \L2V$}}.
  \label{strongsigma}
\Esist
We deduce that $\pi(\phieps)$ and \pier{$p(\phieps)$ converge to $\pi(\phi)$ and $p(\phi)$
strongly in $C^0(\overline Q)$. As $\Ieps(\sigmaeps)= \sigmaeps$ for $\eps\leq 1/\sigmastar$,
we infer that the products $\Ieps(\sigmaeps)\pi(\phieps)$   and $\Ieps(\sigmaeps) p (\phieps)$ 
converge to $\sigma\pi(\phi)$ and $\sigma p(\phi)$} strongly in $\C0H$.
Next, by combining \eqref{strongphi} and \accorpa{convxi}{convzeta}
with well-known monotonicity arguments
(see, e.g., \pier{\cite[Prop.~2.2, p.~38]{Barbu}})
we infer that $\xi\in\beta(\phi)$ and $\zeta\in\sign(\phi-\phistar)$ \aeQ.
Thus, all the equations \accorpa{prima}{quarta} are satisfied.
Finally, \pier{on account of the strong convergence in $\C0H$ for $\phieps$ and~$\sigmaeps$,
it turns out that} $(\phi,\sigma)$ satisfies~\eqref{cauchy} as well.
Therefore, the whole problem \Pbl\ is solved \pier{and Theorem~\ref{Wellposedness} is completely proved.}


\section{Existence of sliding modes}
\label{SLIDINGMODES}
\setcounter{equation}{0}

This section is devoted to the proof of Theorem~\ref{Sliding}.
Our method follows the ideas of \cite{BCGMR}
and relies on a comparison argument for the equation~\eqref{seconda}.
To this end, we define
\Beq
  \xistar := \betaz(\phistar)
  \label{defxistar}
\Eeq
and notice that $\xistar$ is bounded due to~\eqref{hpstar}.
We also remark that $\phiz$, $\phistar$, $\Delta\phistar$ and $\muH$ are bounded as well
and that $\muH$ satisfies $\norma\muH_\infty\leq\norma\muG_\infty$
(see~\eqref{hpz}, \eqref{hpsliding}, \eqref{muGbdd} and~\eqref{maxmuH}).
Furthermore, we introduce some abbreviations
(according to those of Remark~\ref{Calcolo}).
We set for convenience
\Bsist
  && M :=
  \norma\muG_\infty
  + \norma{\Delta\phistar}_\infty
  + \norma\xistar_\infty \,, \quad
  \Mz := \norma{\phiz-\phistar}_\infty 
  \label{defM}
  \\
  && \Mpi := \sup\{|\pi(\phistar(x)+r)| :\ x\in\Omega,\ |r|\leq\Mz\}.
  \label{defMpi}
\Esist
Next, we recall the definition of $\Csys$ given in \eqref{smallness}
and the estimate \eqref{stimamu} for $\mu$ provided by Theorem~\ref{Wellposedness}.
Thus, we have
\Beq
  \norma\mu_\infty
  \leq \Csh \, \frac {2 |\Omega|^{2/3}} \tau \, \rho
  + \hatC 
  = \Csys \, \rho + \hatC .
  \label{mubdd}
\Eeq
Finally, we define $A:(0,+\infty)\to\erre$ by setting
\Beq
  A(\rho) := 
  \Csys \, \rho
  + \hatC
  + M 
  + \Mpi
  \quad \hbox{for $\rho>0$}.
  \label{defArho}
\Eeq
At this point, we set for convenience
\Beq
  \chi := \phi-\phistar
  \label{defchi}
\Eeq
and write the equation \pier{\eqref{seconda} in the form}
\Beq
  \tau \dt\chi
  - \Delta\chi
  + \xi - \xistar
  + \rho \, \zeta
  = \mu + \muH \pier{{}+\Delta \phistar} - \xistar - \pi(\phistar+\chi).
  \label{secondasl}
\Eeq
By accounting for the above definitions of $M$ and $\Mpi$ and using~\eqref{mubdd},
we obtain the inequalities (which we use later~on)
\Bsist
  & \tau \dt\chi
  - \Delta\chi
  + \xi - \xistar
  + \rho \, \zeta
  \hskip -.5em
  & \leq \Arho - \Mpi - \pi(\phistar+\chi)
  \quad \aeQ
  \label{secondapos}
  \\
  & \tau \dt\chi
  - \Delta\chi
  + \xi - \xistar
  + \rho \, \zeta
  \hskip -.5em
  & \geq - \Arho + \Mpi - \pi(\phistar+\chi)
  \quad \aeQ .
  \label{secondaneg}
\Esist
At this point, we recall that $\Csys<1$ by \eqref{smallness}
and define $\rhostar$ \pier{in this way}
\Beq
  \rhostar :=
  \frac 1 {1-\Csys} \Bigl(
    \hatC 
    + M
    + \Mpi
    + \frac \tau T \, \Mz  
  \Bigr) 
  \label{defrhostar}
\Eeq
so that 
\Beq
  \Arho + \frac \tau T \, \Mz 
  < \rho
  \quad \hbox{for every $\rho>\rhostar$}.
  \label{disugArho}
\Eeq
Once such a definition is given, we fix $\rho>\rhostar$
and a solution $\soluz$ to problem \Pbl\ 
whose component $\mu$ satisfies~\eqref{stimamu}, i.e., \eqref{mubdd},
and prove the property of the statement.
To this end, we notice that \eqref{disugArho} implies $\Arho<\rho$.
Thus, the definition
\Beq
  \Tstar := \frac \tau {\rho-\Arho} \, \Mz
  \label{defTstar}
\Eeq
makes sense and yields $\Tstar\in[0,T)$.
At this point, we prove that $\chi(t)=0$, i.e., $\phi(t)=\phistar$, 
for every $t\in[\Tstar,T]$ by taking advantage of a comparison argument.
We consider the auxiliary problem of finding a pair $(w,\eta)$ satisfying
\Bsist
  && w \in H^1(0,T), \quad  
  \eta \in L^\infty(0,T)
  \aand
  \eta \in \sign w
  \quad \aet
  \label{regw}
  \\
  && \tau w' + \rho \, \eta = \Arho
  \quad \aet
  \aand
  w(0) = \Mz \,.
  \label{odew}
\Esist
Since $\sign$ is a maximal monotone operator, 
the above problem has a unique solution $(w,\eta)$.
Moreover, by noticing that $\Arho/\rho\in[0,1)\subset\sign 0$,
we see by a simple computation that 
\Beq
  w(t) = \Bigl( \Mz - \frac {\rho-\Arho} \tau \, t \Bigr)^+
  \quad \hbox{for every $t\in[0,T]$}
  \label{formulaw}
\Eeq
(with $\eta(t)=1$ and $\eta(t)=\Arho/\rho$ for $t<\Tstar$ and $t>\Tstar$, respectively).
As $w(t)=0$ for every $t\in[\Tstar,T]$,
it suffices to prove that $|\chi(x,t)|\leq w(t)$ for every $(x,t)\in Q$.
To this end, we read $w$ and $\eta$ as space independent functions defined in~$Q$
and the equation in \eqref{odew}~as
\Beq
  \tau \dt w - \Delta w + \rho \, \eta = \Arho
  \quad \aeQ
  \label{pblw}
\Eeq
with the boundary condition $\dn w=0$ on~$\Sigma$.
Since $0\leq w\leq\Mz$ by~\eqref{formulaw},
we see that $|\pi(\phistar\pm w)|\leq\Mpi$ by the definition \eqref{defMpi} of~$\Mpi$.
Therefore, $w$~also satisfies the inequalities
\Bsist
  && \tau \dt w - \Delta w + \rho \, \eta
  \geq \Arho - \Mpi - \pi(\phistar+w)
  \quad \aeQ
  \label{pblpos}
  \\
  && \tau \dt w - \Delta w + \rho \, \eta
  \geq \Arho - \Mpi + \pi(\phistar-w)
  \quad \aeQ .
  \label{pblneg}
\Esist
At this point, we start with the announced comparison argument.
We recall the inequalities \eqref{secondapos} and \eqref{pblpos}
and take the difference.
We have
\Beq
  \tau \dt(\chi-w)
  - \Delta(\chi-w)
  + \xi - \xistar
  + \rho \, (\zeta-\eta)
  \leq \pi(\phistar+w) - \pi(\phistar+\chi)
  \non
\Eeq
and multiply this inequality by the positive part $(\chi-w)^+$.
By also integrating over~$Q_t$
and noticing that $(\chi-w)^+(0)=0$ by the definition \eqref{defM} of~$\Mz$
and the initial condition $w(0)=\Mz$, we obtain
\Bsist
  &&
  \pier{\frac \tau 2} \iO |(\chi-w)^+(t)|^2 
  + \intQt |\nabla(\chi-w)^+|^2
  \non
  \\
  && \quad {}
  + \intQt (\xi-\xistar) (\chi-w)^+
  + \rho \intQt (\zeta-\eta) (\chi-w)^+
  \non
  \\
  && \leq \Lpi \intQt |(\chi-w)^+|^2 .
  \non
\Esist
Next, we show that the third and \pier{fourth} term on the \lhs\ 
of the above equality are nonnegative by using the monotonicity of $\beta$ and~$\sign$.
Indeed, in the set $Q^+$ where $(\chi-w)^+\not=0$ we have $\chi>w$,
whence $\zeta\geq\eta$ since $\zeta\in\sign(\chi)$ and $\eta\in\sign(w)$.
In the same set, we also have
$\phi=\phistar+\chi>\phistar+w\geq\phistar$,
so that $\xi\geq\xistar$ since $\xi\in\beta(\phi)$ and $\xistar\in\beta(\phistar)$.
Therefore, we can apply the Gronwall lemma and conclude that $(\chi-w)^+=0$ \aeQ, \pier{i.e., $\chi\leq w$}.
Now, we consider the sum of the inequalities \eqref{secondaneg} and~\eqref{pblneg}, \pier{that~is},
\Beq
  \tau \dt(\chi+w)
  - \Delta(\chi+w)
  + \xi - \xistar
  + \rho \, (\zeta+\eta)
  \geq \pi(\phistar-w) - \pi(\phistar+\chi)
  \non
\Eeq
and multiply it by $-(\chi+w)^-$.
We analogously have that
\Bsist
  &&
 \pier{\frac \tau 2} \iO |(\chi+w)^-(t)|^2 
  + \intQt |\nabla(\chi+w)^-|^2
  \non
  \\
  && \quad {}
  - \intQt (\xi-\xistar) (\chi+w)^-
  - \rho \intQt (\zeta+\eta) (\chi+w)^-
  \non
  \\
  && \leq \Lpi \intQt |(\chi+w)^-|^2 
  \non
\Esist
and we can treat the integrals involving the nonlinearities similarly as before.
In the set $Q^-$ where $(\chi+w)^-\not=0$ we have $\chi<-w$,
whence $\zeta\leq-\eta$ since $\zeta\in\sign(\chi)$ and $-\eta\in\sign(-w)$.
In the same set, we also have
$\phi=\phistar+\chi<\phistar-w\leq\phistar$,
so that $\xi\leq\xistar$
since $\xi\in\beta(\phi)$ and $\xistar\in\beta(\phistar)$.
Therefore, we can apply the Gronwall lemma also in this case
and deduce that $(\chi+w)^-=0$ \aeQ, i.e., that $-\chi\leq w$.
By combining the inequalities we have obtained, we conclude that
$|\chi|\leq w$, i.e., that $|\phi-\phistar|\leq w$.
\pier{Therefore,} this and \eqref{formulaw} imply that $\phi(t)=\phistar$ for $t\geq\Tstar$,
and the proof is complete.



\section*{Acknowledgements}
\pier{This research activity has been performed in the framework of an
Italian-Romanian  {three-year project on ``Control and 
stabilization problems for phase field and biological systems'' financed by the Italian CNR and the Romanian Academy.} 
Moreover, the financial support of  the project Fondazione Cariplo-Regione Lombardia  MEGAsTAR 
``Matema\-tica d'Eccellenza in biologia ed ingegneria come \bettibis{acceleratore} 
di una nuova strateGia per l'ATtRattivit\`a dell'ateneo pavese'' is gratefully acknowledged by the authors. 
The present paper 
also benefits from the support of the MIUR-PRIN Grant 2015PA5MP7 ``Calculus of Variations'' for PC and GG, 
the GNAMPA (Gruppo Nazionale per l'Analisi Matematica, la Probabilit\`a e le loro Applicazioni) 
of INdAM (Istituto Nazionale di Alta Matematica) for PC, GG and~ER,
and the UEFISCDI project PN-III-ID-PCE-2016-0011 for GM.}



\vspace{3truemm}

{\small%

\Begin{thebibliography}{99}

\bibitem{BFPU08}
G. Bartolini, L. Fridman, A. Pisano, E. Usai (eds.),
``Modern Sliding Mode Control \revis{Theory New} Perspectives and Applications'',
Lecture Notes in Control and Information Sciences {\bf 375},
Springer, 2008.

\bibitem{Barbu}
V. Barbu,
``Nonlinear \revis{Differential Equations of Monotone Types in Banach Spaces}'',
Springer, New York, 2010.

\bettibis{\bibitem{BCGMR}
V. Barbu, P. Colli, G. Gilardi, G. Marinoschi, E. Rocca,
Sliding mode control for a nonlinear phase-field system,
{\it SIAM J. Control Optim.} {\bf 55} (2017), 2108--2133.}

\bibitem{BCT}
E. Bonetti, P. Colli and  G. Tomassetti, 
A non-smooth regularization of a forward-backward parabolic equation,
{\it Math. Models Methods Appl. Sci.} {\bf 27} (2017), 641--661.

\bibitem{Brezis}
H. Brezis,
``Op\'erateurs \revis{Maximaux Monotones et Semi-groupes de Contractions
dans les Espaces} de Hilbert'',
North-Holland Math. Stud.
{\bf 5},
North-Holland,
Amsterdam,
1973.

\bibitem{BG87}
F.\ {B}rezzi, G.\ {G}ilardi,
\newblock {Chapters 1-3 in ``Finite \revis{Element Handbook}'',
\newblock H.\ {K}ardestuncer and D.\ H.\ {N}orrie (Eds.)},
McGraw--Hill Book Co., New York, 1987.

\bibitem{BC} H. M. Byrne, M. A. J. Chaplain, 
Growth of nonnecrotic tumors in the presence and absence of inhibitors,
{\it Math. Biosci.} {\bf 130} (1995), 151--181.  

\bibitem{CRS11}
\revis{M.-B.} Cheng, V. Radisavljevic, W.-C. Su,
Sliding mode boundary control of a parabolic PDE system with parameter variations and boundary uncertainties, 
{\it Automatica J. IFAC} {\bf 47} (2011), 381--387.

\bibitem{CGH}
P.~Colli, G.~Gilardi, D.~Hilhorst, On a Cahn--Hilliard type phase field model related to tumor growth, 
{\it Discrete Contin. Dyn. Syst.} {\bf 35} (2015), 2423--2442.

\bibitem{CGRS1}
P. Colli, G. Gilardi, E. Rocca, J. Sprekels, 
Vanishing viscosities and error estimate for a Cahn--Hilliard type phase-field system related to tumor growth, 
{\it Nonlinear Anal. Real World Appl.} {\bf 26} (2015), 93--108.

\bibitem{CGRS2}
P.~Colli, G.~Gilardi, E.~Rocca, J.~Sprekels, 
Asymptotic analyses and error estimates for a Cahn--Hilliard type phase field system modelling tumor growth, {\it Discrete Contin. Dyn. Syst. Ser. S.} {\bf 10} (2017), 37--54.

\bibitem{CGRS17}
P.~Colli, G.~Gilardi, E.~Rocca, J.~Sprekels, 
Optimal distributed control of a diffuse interface model of tumor growth, 
{\it Nonlinearity} {\bf 30} (2017), 2518--2546. 

\bibitem{Colt1}
M.~Colturato, 	
Solvability of a class of phase field systems related to a sliding mode control problem, 
{\it Appl. Math.} \textbf{6} (2016), 623--650. 
\bibitem{Colt2} 
M.~Colturato, 
On a class of conserved phase field systems with a maximal monotone perturbation, 
{\it Appl. Math. Optim.} \revis{\textbf{4} (2017), 1--35.}

\bibitem{CL10} V. Cristini, J. Lowengrub, 
``Multiscale \revis{Modeling of Cancer}. An Integrated Experimental and Mathematical Modeling Approach'', Cambridge Univ. Press, 2010. 

\bibitem{CLLW09} V. Cristini, X. Li, J.S. Lowengrub, S.M. Wise, 
Nonlinear simulations of solid
tumor growth using a mixture model: \revis{invasion} and branching,
{\it J. Math. Biol.} {\bf 58} (2009), 723--763. 

\bibitem{DFRSS17}
M.~Dai, E.~Feireisl, E.~Rocca, G.~Schimperna, M.~Schonbek, 
Analysis of a diffuse interface model for multispecies tumor growth, 
{\it Nonlinearity} {\bf 30} (2017), 1639--1658.

\bibitem{EFF06}
C. Edwards, E. Fossas Colet, L. Fridman (eds.),
``Advances in Variable Structure and Sliding Mode Control'',
Lecture Notes in Control and Information Sciences {\bf 334},
Springer-Verlag, 2006.

\bibitem{ES99}
C. Edwards, S. Spurgeon,
``Sliding Mode Control: Theory and Applications'',
Taylor and Francis, London, 1999.

\bibitem{FMI11}
L. Fridman, J. Moreno, R. Iriarte (eds.),
``Sliding Modes After the First Decade of the 21st Century: State of the Art'',
Lecture Notes in Control and Information Sciences {\bf 412},
Springer, 2011.

\bibitem{FGR}
S.~Frigeri, M.~Grasselli, E.~Rocca, On a diffuse interface model of tumor growth, 
{\it European J. Appl. Math.} {\bf 26} (2015), 215--243.

\bibitem{FLR}
S. Frigeri, K.F. Lam, E. Rocca, 
On a diffuse interface model for tumour growth 
with non-local interactions and degenerate mobilities, 
\revis{to appear in} ``Solvability,
Regularity, Optimal Control of Boundary Value Problems for PDEs'',
P.~Colli, A.~Favini, E.~Rocca, G.~Schimperna, J.~Sprekels~(ed.),
Springer INdAM Series, Springer, Milan, 2017 \revis{(see also
preprint arXiv:1703.03553 (2017), pp.~1--28)}

\bibitem{GLDarcy}
H.~Garcke, K.F.~Lam, Global weak solutions and asymptotic limits of a Cahn--Hilliard--Darcy system modelling tumour growth, {\it AIMS Mathematics} { \bf 1} (2016), 318--360.

\bibitem{GLDirichlet}
H.~Garcke, K.F.~Lam, Analysis of a Cahn--Hilliard system with non zero Dirichlet conditions modelling tumour growth with chemotaxis,  {\it Discrete Contin. Dyn. Syst.} {\bf 37} (2017), 4277--4308.

\bibitem{GLNeumann}
H.~Garcke, K.F.~Lam, 
Well-posedness of a Cahn--Hilliard system modelling tumour growth with chemotaxis and active transport, {\it European J. Appl. Math.} {\bf 28} (2017), 284--316.

\bibitem{GLNS} 
H. Garcke, K.F. Lam, \revis{R. N\"urnberg}, E. Sitka, 
A multiphase Cahn--Hilliard--Darcy model for tumour growth with necrosis, 
preprint arXiv:1701.06656v1 [math.AP] (2017), \revis{pp.}~1--43. 

\bibitem{GLR17}
H.~Garcke, K.F.~Lam, E.~Rocca, 
Optimal control of treatment time in a diffuse interface model for tumor growth,
{\it Appl. Math. Optim.}, DOI: 10.1007/s00245-017-9414-4 (2017). 

\bibitem{GLSS16}
H. Garcke, K.F. Lam, E. Sitka, V. Styles, 
A \revis{Cahn--Hilliard--Darcy} model for
tumour growth with chemotaxis and active transport, 
{\it Math. Models Methods Appl. Sci.} {\bf 26} (2016),
1095--1148. 

\bibitem{HZO12}
A. Hawkins-Daarud, K.G. van der Zee, J.T. Oden, 
Numerical simulation of a thermodynamically
consistent four-species tumor growth model, 
{\it Int. J. Numer. Methods Biomed. Eng.}  
{\bf 28} (2012), 3--24. 

\bibitem{HKNZ15} 
D. Hilhorst, J. Kampmann, T.N. Nguyen, K.G. van der Zee, 
Formal asymptotic limit of a diffuse-interface tumor-growth model, 
{\it Math. Models Methods Appl. Sci.} {\bf  25} (2015), 1011--1043. 

\bibitem{I76}
U. Itkis,
``Control \revis{Systems of Variable Structure}'',
Wiley, 1976.

\revis{\bibitem{LO02}
L. Levaggi,
Infinite dimensional systems' sliding motions,
{\it \pier{Eur. J. Control}}~{\bf 8} (2002), 508--516.}

\bibitem{Levaggi13}
L. Levaggi, 
Existence of sliding motions for nonlinear evolution equations in Banach spaces, 
\revis{{\it Discrete Contin. Dyn. Syst.} 2013, ``Dynamical Systems, Differential Equations and Applications. 9th AIMS Conference. Suppl.'', pp.~477--487.}

\bibitem{MR}
S.~Melchionna, E.~Rocca, 
Varifold solutions of a sharp interface limit of a diffuse interface model for tumor growth, {\it Interfaces and Free Bound.} to appear (2017) (see also
preprint arXiv:1610.04478 [math.AP] (2016), pp.~1--29)

\bibitem{O83}
Y.V. Orlov,
Application of Lyapunov method in distributed systems,
{\it Autom. Remote Control} {\bf 44} (1983), 426--430.

\bibitem{O00}
Y.V. Orlov, 
Discontinuous unit feedback control of uncertain infinite dimensional
systems, 
{\it IEEE Trans. Automatic Control} {\bf 45} (2000), 834--843.

\bibitem{OU83}
Y.V. Orlov, V.I.  Utkin,
 Use of sliding modes in distributed system control
problems,
{\it Autom. Remote Control} {\bf 43} (1983), 1127--1135.

\bibitem{OU87}
Y.V. Orlov, V.I.  Utkin,
Sliding mode control in indefinite-dimensional systems, 
{\it Automatica
J. IFAC} {\bf 23} (1987), 753--757.

\bibitem{OU98}
Y.V. Orlov, V.I.  Utkin, 
{Unit sliding mode control in infinite-dimensional systems. Adaptive
learning and control using sliding modes},  
{\it Appl. Math. Comput. Sci.} {\bf 8} (1998), 7--20.

\bibitem{RS}
E.~Rocca, R.~Scala,
A rigorous sharp interface limit of a diffuse interface model related to tumor growth,
{\it J. Nonlinear Sci.} {\bf 27} (2017), 847--872. 

\bibitem{Simon}
J. Simon,
{Compact sets in the space $L^p(0,T; B)$},
{\it Ann. Mat. Pura Appl.~(4)\/}
{\bf 146} (1987), 65--96.

\bibitem{Utkin92}
V. Utkin,
``Sliding Modes in Control and Optimization'',
Communications and Control Engineering Series,
Springer-Verlag, 1992.

\bibitem{UGS09}
V. Utkin, J. Guldner, J. Shi,
``Sliding Mode Control in Electro-Mechanical Systems'', 2nd Edition
CRC Press, Automation and Control Engineering, 2009.

\bibitem{XLGK13}
H. Xing, D. Li, C. Gao, Y. Kao,
Delay-independent sliding mode control for a class of quasi-linear parabolic distributed parameter systems with time-varying delay,
{\it J. Franklin Inst.} {\bf 350} (2013), 397--418.

\bibitem{WLFC08} S.M. Wise, J.S. Lowengrub, H.B. Frieboes, V. Cristini, 
Three-dimensional multispecies nonlinear tumor growth--I: 
model and numerical method, 
{\it J. Theoret. Biol.} {\bf 253} (2008),  524--543.

\bibitem{YO99}
K.D. Young, \"U. \"Ozg\"uner (eds.),
``Variable Structure Systems, Sliding Mode and Nonlinear \revis{Control''},
Springer-Verlag, 1999.

\End{thebibliography}

}

\End{document}

\bye